\documentclass[12pt]{amsart}
\usepackage{amsmath,mathrsfs}
\usepackage{amssymb}
\usepackage{latexsym}
\usepackage{amsfonts}
\usepackage{hyperref}
\usepackage{graphicx}
\usepackage{psfrag}
\usepackage{pgfplots}
\usepackage{verbatim}
\usepackage{tikz}
\usetikzlibrary{arrows,calc,chains,shapes,dsp}
\usepackage[siunitx]{circuitikz}
\usepackage{psfrag}
\usepgfplotslibrary{fillbetween}
\def\C{\mathbb C}

\def\R{{\mathbb R}}
\def\neq{\not=}

\def\@tempb{saamsart}
\renewcommand{\theequation}{\thesection.\arabic{equation}}
\newtheorem{Pa}{Paper}[section]
\newtheorem{Tm}[Pa]{{\bf Theorem}}
\newtheorem{La}[Pa]{{\bf Lemma}}
\newtheorem{Ob}[Pa]{{\bf Observation}}
\newtheorem{Cy}[Pa]{{\bf Corollary}}

\newtheorem{Pn}[Pa]{{\bf Proposition}}

\newtheorem{Ex}[Pa]{{\bf Example}}
\newtheorem{Dn}[Pa]{{\bf Definition}}
\author[D. Alpay]{Daniel Alpay}
\address{(DA) Schmid College of Science and Technology\newline
Chapman University\newline
One University Drive
Orange, California 92866\\
USA}
\email{alpay@chapman.edu}
\author[I. Lewkowicz]{Izchak Lewkowicz}
\address{(IL) Department of Electrical Engineering
\newline
Ben Gurion University of the Negev \newline P.O.B. 653,
\newline
Be'er Sheva 84105, \newline ISRAEL }
\email{izchak@bgu.ac.il}
\date{}
\title[``Wrong" side interpolation]
{``Wrong" side interpolation by low degree positive real
rational functions}

\begin{document}
\bibliographystyle{plain}
\begin{abstract}
Using polynomial interpolation, along with structural properties of
the family of rational positive real functions, we here show that a
set of $m$ nodes in the open {\em left} half of the complex plane,
can always be mapped to anywhere in the complex plane by rational
positive real functions whose degree is at most $m$. Moreover we
introduce an ~{\em easy-to-find}~ parametrization
in $\R^{2m+3}$ of a large subset of these interpolating functions.
\end{abstract}
\maketitle
\renewcommand{\theequation}
{\thesection.\arabic{equation}}

\section{Introduction}
\setcounter{equation}{0}

{\bf Problem Formulation}
\vskip 0.2cm

A framework for many classical interpolation problems is as
follows. Given a set of distinct nodes $x_1,~\ldots~,~x_m$,
image points $y_1,~\ldots~,~y_m$ (not necessarily distinct)
and a family of functions $\mathcal{F}$, find whether there
exist functions $f\in\mathcal{F}$ so that
\begin{equation}\label{eq:interp}
y_j=f(x_j)\quad\quad\quad\quad j=1,~\ldots~,~m.
\end{equation}
If yes, parameterize all of them, preferably within a degree
bound. There is a vast literature on the subject see e.g.
\cite{AA}-\cite{BoDy}, \cite{Ko}-\cite{YS}.
To simplify the discussion, we here focus on scalar real rational
functions. Thus degree simply means the maximum between the degree
of the numerator and of the denominator polynomials. The
polynomial (a.k.a. the Lagrange) interpolation \cite{La} (in
\cite{Mej} it is attributed to \cite{Wa}) is probably the best
known problem in this framework. For the case where $\mathcal{F}$
is the set of rational functions see \cite{ABKW} and if in addition
all functions in $\mathcal{F}$ are analytic in a disk of a
prescribed radius in $\C$, the problem was addressed in \cite{AA}.
\vskip 0.2cm

We shall denote by $\C_r$ $(\overline{\C}_r)$ the open (closed)
right half plane (the subscript stands for ``right"). The family of
functions $\mathcal{F}$ we here focus on, is of ~{\em positive
real},~ i.e. analytically mapping the open right half plane to its
closure. Namely, a real rational function $~f(s)~$ of a complex
variable $~s~$ is said to be positive if
\begin{equation}\label{eq:DefPos}
{\rm Re}\left(f(s)\right)\geq 0
\quad\quad\quad\quad\forall s\in\C_r~.
\end{equation}
Interpolation problem with rational positive real functions
can be further classified by the domain the nodes
$x_1,~\ldots~,~x_m$ belong to.
\vskip 0.2cm

If the nodes $x_j$ are in $\C_r$, this amounts to the classical
Nevanlinna-Pick interpolation problem, see e.g.
\cite[Theorem 18.1]{BGR} and for real functions \cite{YS}. There,
from the interpolation data one constructs the
Pick matrix whose $j,k$ element is given by
\[
\frac{y_j^*+y_k}{x_j^*+x_k}
\quad\quad\quad\quad j,~k=1,~\ldots~,~m.
\]
It is known that there exist interpolating functions if and only
if the Pick matrix is positive semi-definite. Moreover, all
interpolating functions may be parameterized through this Pick
matrix. Recall that having the Pick matrix positive semi-definite
implies that each 2-dimensional minor is non-negative, which in turn
can be written as,
\begin{equation}\label{eq:NevPick}
\frac{|x_j-x_k|^2}{{\rm Re}(x_j){\rm Re}(x_k)}
\geq
\frac{|y_j-y_k|^2}{{\rm Re}(y_j){\rm Re}(y_k)}
\quad\quad\quad\quad m\geq j>k\geq 1.
\end{equation}
This condition means that the map from the nodes $x_1~,~\ldots~,~x_m$ to
the image points $y_1~,~\ldots~,~y_m$, is contractive in $\C_r$ in
the sense of Eq. \eqref{eq:NevPick}. This illustrates the fact that
the interpolation problem cannot be solvable for arbitray data set.
\vskip 0.2cm

\noindent
If the nodes are confined to the imaginary axis, an interpolation scheme,
elegant in its simplicity, appeared in \cite{ZehebLempel1966}.
\vskip 0.2cm

\noindent
If the nodes $x_j$ are in $\overline{\C}_r$ (with possibly some nodes
on $i\R$) the problem is much harder, see e.g. \cite[Chapter 21]{BGR},
\cite{BoCa}, \cite{BoDy}, \cite{Ko} and \cite{Sa}.
\vskip 0.2cm

\noindent
If the interpolation data is in whole plane, provided that
\begin{equation}\label{eq:OddNevPick}
{\rm Re}(x_j){\rm Re}(y_j)
>0\quad\quad\quad j=1,~\ldots~,~m,
\end{equation}
one can still resort to the
classical Nevanlinna-Pick interpolation scheme: First, complete
the data set so that if $x, y$ is an interpolation pair,
then so is $-x, -y$. Then, from this extended data, take the $m$
nodes which are in $\C_r$, construct the the corresponding
Pick matrix and proceed as usual. Finally, use the fact (see
\cite{YS}) that whenever the Pick matrix is positive
semi-definite, among the interpolation functions there exists
some function with odd symmetry, i.e. $f(s)=-f(-s)~$ (a.k.a. Foster or
lossless functions, see e.g. \cite{Belev1968}, \cite{CL}, \cite{YS}).
If instead of the left and right half planes, $\C$ is partitioned
to the unit disk and its exterior, a similar idea is presented
in \cite[Section 5]{AA}.
\vskip 0.2cm

{\em In this work we focus on the case where the nodes
$x_1,~\ldots~,~x_m$ are all in $\C_l$ (the open left half
plane). We parameterize a large subset of rational positive real
interpolating functions whose degree is less or equal to
$m$. In particular, it is shown that this set is never empty.}
\vskip 0.2cm

\noindent
A key idea is the following: We construct two rational functions
sharing the same denominator: (i) an interpolating
function $p(s)$ (not necessarily positive real) and (ii) a strictly positive
real rational function ${\scriptstyle\Delta}(s)$ vanishing at the nodes.
Thus, for all $r\in\R$ the parametric rational function
\begin{equation}\label{eq:KeyIdea}
f(s)=p(s)+r{\scriptstyle\Delta}(s),
\end{equation}
is interpolating. Moreover, for ${ r}$
``sufficiently large", $f(s)$ turns to be positive real.
\vskip 0.2cm

\noindent
Interestingly, we can mention two ideas conceptually similar to those
in the current work, which have appeared within completely different
interpolation frameworks: (i) The fact that for interpolation by low
degree rational functions, one should separately treat numerators and
denominators, appeared in the context of Schur functions in
\cite[Theorems 1, 2]{Georgiou1999}. (ii) In
\cite{AlpayLew2014} we have used Eq. \eqref{eq:KeyIdea}
for interpolation by structured matrix-valued polynomials.
For example, where $p(s)$ was an interpolating polynomial which on
$i\R$, attained Hermitian values, while the polynomial
${\scriptstyle\Delta}(s)$, vanishing at the nodes, was positive definite
on $i\R$, see \cite[Eq. (1.6)]{AlpayLew2014}.
\vskip 0.2cm

\noindent
In Section II we present a five steps interpolation procedure:
\vskip 0.2cm

\noindent
In Step 1 we parameterize all candidates for denominator polynomials
of the sought interpolating functions. Namely, all real polynomials,
of degree of at most $~m$, non-vanishing at the nodes.
\vskip 0.2cm

\noindent
In Step 2,
to each of these denominator polynomials we match a numerator to
obtain $p(s)$, a rational interpolating function (not
necessarily positive real).
\vskip 0.2cm

\noindent
In Step 3, we construct ${\scriptstyle\Delta}(s)$ strictly
positive real rational functions, vanishing at the nodes. We now
restrict the denominators of $p(s)$, the rational functions from
Step 2 to the subset of the resulting deniminators of
${\scriptstyle\Delta}(s)$.
\vskip 0.2cm

\noindent
To each of the resulting interpolating function $p(s)$, we
add $r{\scriptstyle\Delta}(s)$, a weighted version of the
strictly positive real rational functions, vanishing at the nodes
(which shares the same denominar). Thus,
$p(s)+ r{\scriptstyle\Delta}(s)$ is an interpolating function,
of degree of at most $m$. Furthermore,
for $r$ ``sufficiently large" it is positive real.
\vskip 0.2cm

\noindent
A closer scrutiny reveals that all interpolating, positive real
rational functions obtained, are so that the degree of the
denominator is larger or equal to the degree of the numerator.
\vskip 0.2cm

\noindent
In Step 5, we complete our the description of positive real
interpolating functions as follows: We repeat the previous
steps by constructing positive real interpolating functions from
the original nodes $~x_j$ {\em but to}\begin{footnote}{Consider
the case where $y_j=0$ for some $j$.}\end{footnote}~ $\frac{1}{y_j}$.
Finally, as the sought solution, we take the reciprocal of these
functions.
\vskip 0.2cm

\noindent
In Section III we illustrate the above procedure by detailed examples.
and add concluding remarks.

\section{A Recipe}
\label{sec:Recipe}

\subsection{Step 1: Constructing all real monic polynomials,
of degree ${\mathbf m}$ and ${\mathbf m-1}$, non-vanishing
at the nodes}
\label{SubSec:Denom}
\setcounter{equation}{0}
\vskip 0.2cm

\noindent
We shall do it in stages.
\vskip 0.2cm

\noindent
{\bf 1a}\quad Constructing all complex polynomials of degree of at
most~ $m$ ~with no roots at prescribed distinct points~
\mbox{$x_1,~\ldots~,~x_m\in\C$}.
\vskip 0.2cm

\noindent
First, denote by~ $\eta(s)$ ~the monic polynomial
(of degree~ $m$) ~whose roots are the prescribed distinct points~
\mbox{$x_1,~\ldots~,~x_m\in\C$},
\begin{equation}\label{eq:eta}
\eta(s):=\prod\limits_{j=1}^m(s-x_j).
\end{equation}
Next, denote by~ \mbox{$\phi_1(s),~\ldots~,~\phi_m(s)$} ~the monic
divisors of $~\eta(s)~$ of degree $~m-1$, i.e.
\begin{equation}\label{eq:phi}
\phi_j(s)=
{\scriptstyle
\frac{\eta(s)}{s-x_j}
}
=
\prod\limits_{\begin{smallmatrix}k=1\\k\not=j
\end{smallmatrix}}^m(s-x_k)\quad\quad\quad
j=1,~\ldots~,~m.
\end{equation}
\begin{La}\label{LemmaTilde}
For distinct $x_1,~\ldots~,~x_m\in\C~$ let
$\eta(s)$ and $\phi_1(x)~,~\ldots~,~\phi_m(x)$
be as in \eqref{eq:eta} and \eqref{eq:phi}, respectively.
\vskip 0.2cm

1. The set of all polynomials of degree of
at most $m-1$ can be parameterized by
\[
\sum\limits_{j=1}^mc_j\phi_j(s)
\quad\quad\quad\quad c_j\in\C.
\]
2. The set of all polynomials of degree of at
most $m$ can be parameterized by
\[
\tilde{d}(s)=b\eta(s)+\sum\limits_{j=1}^mc_j{\phi}_j(s)
\quad\quad\quad\quad b, c_j\in\C.
\]
3. If in addition,
\[
c_1{\cdots}c_m\neq 0,
\]
$\tilde{d}(s)$ forms the set of all polynomials of degree of at
most $m$ which do not vanish at \mbox{$x_1~,~\ldots~,~x_m$}.
\end{La}

{\bf Proof}\quad 1.
For each $j$, $j\in[1,~m]$, let
\[
{\phi}_j(s)=\sum\limits_{k=0}^{m-1}a_{jk}s^k
\]
be identified with
\[
a_j:=
\left(\begin{matrix}a_{jo}\\ a_{j1}\\ \vdots\\ a_{j,m-1}
\end{matrix}\right).
\]
Now by construction,
\[
\begin{split}
\begin{pmatrix}
1&x_1&x_1^2&\ldots&x_1^{m-1}\\
1     &x_2   &x_2^2 &\ldots&x_2^{m-1}\\
\vdots&\vdots&\vdots&\vdots&\vdots   \\
1     &x_m   &x_m^2 &\ldots&x_m^{m-1}
\end{pmatrix}
\begin{pmatrix}
a_{10}&a_{20}&\ldots&a_{m0}\\
a_{11}&a_{21}&\ldots&a_{m1}\\
a_{12}&a_{22}&\ldots&a_{m2}\\
\vdots&\vdots&\vdots&\vdots\\
a_{1,m-1}&a_{2,m-1}&\ldots&a_{m,m-1}
\end{pmatrix}&=\\
&\hspace{-3cm}
=\begin{pmatrix}
*&0&0&\ldots&0\\
0&*&0&\ldots&0\\
0&0&*&\ldots&0\\
\vdots&\vdots&\vdots&\vdots&\vdots\\
0&0&0&\ldots&*
\end{pmatrix}
\end{split}
\]
where $~*~$ denotes a non-zero element.
\vskip 0.2cm

The leftmost is a Vandermonde matrix, which is non-singular,
whenever \mbox{$x_1~,~\ldots~,~x_m$} are distinct, see e.g.
\cite[Item 0.9.11]{HJ}.
Thus, the middle matrix whose entries are $a_{jk}$ must
be non-singular, so in particular its columns form
a basis to $\C^m$. In other words
\mbox{${\phi}_1(s),~\ldots~,~{\phi}_m(s)$} form
a basis to all polynomials of degree of at most $m-1$.
\vskip 0.2cm

2. In a similar way, the coefficients of $\eta(s)$ in \eqref{eq:eta}
can be identified with a vector in $\C^{m+1}$. By adding a
bottom zero to each the above vectors $a_1~,~\ldots~,~a_m$,
they are embedded in $\C^{m+1}$. Thus, the problem is reduced
to verifying that the last element in the vector associated with
$\eta(s)$ is non-zero, but by  construction $\eta(s)$ is of
degree $m$, so indeed its bottom element is non-zero.
\vskip 0.2cm

3. Having no roots at \mbox{$x_1~,~\ldots~,~x_m$}. Note
that,
\[
{\tilde{d}}(s)_{|_{s=x_j}}=c_j{\phi}_j(s)_{|_{s=x_j}}=
\left(c_j\prod\limits_{\begin{smallmatrix}k=1\\k\not=j
\end{smallmatrix}}^m(x_j-x_k)\right)_{|_{c_j\not=0}}\not=0.
\]
In other words, having $~c_1{\cdots}c_m\neq 0~$ is
necessary and sufficient for these polynomials not to vanish at
the original points \mbox{$x_1~,~\ldots~,~x_m$}.
Thus the claim is established.
\qed
\vskip 0.2cm

In the sequel we focus on real polynomials and real rational functions.
Hence, we have the next stage.
\vskip 0.2cm

\noindent
{\bf 1b}\quad Assuming the prescribed distinct points
\mbox{$x_1,~\ldots~,~x_m\in\C$}, are closed under complex
conjugation, constructing all real polynomials of degree of at
most~ $m$, ~with no roots at these points.
\vskip 0.2cm

\noindent
Assume hereafter that, if necessary, the original set of
distinct points \mbox{$x_1,~\ldots~,~x_m\in\C$} ~is complemented
so it is closed under complex conjugation. Namely,
\begin{equation}\label{eq:ComplexConjugation}
{\rm Im}(x_j)>0~\Longrightarrow~ x_{j+1}=x_j^*.
\end{equation}
Note that then in Eq. \eqref{eq:eta} the resulting~ $\eta(s)$ ~is
real.
\vskip 0.2cm

We now construct the sought polynomials.

\begin{La}\label{La:Denom}
For distinct~ $x_1,~\ldots~,~ x_m\in\C$, ~closed under complex
conjugation, let~ $\eta(s)$ ~and~
\mbox{$\phi_1(s),~\ldots~,~\phi_m(s)$} ~be as in Eqs.
\eqref{eq:eta} and \eqref{eq:phi}, respectively.
The set of all real polynomials $~\tilde{d}(s)~$ of degree, of at most
$~m$, with no roots at the
original points \mbox{$x_1,~\ldots~,~x_m$}, can be parametrized by
\begin{equation}\label{eq:TildeD}
\tilde{d}(s)=b\eta(s)+\sum\limits_{j=1}^mc_j{\phi}_j(s)
\quad\quad
{\scriptstyle b}\in\R
\quad\quad
\begin{smallmatrix}
{\rm Im}(x_j)>0&~&0\not=&c_{j+1}=c_j^*
\\~\\
{\rm Im}(x_j)=0&~&0\not=&c_j\in\R.
\end{smallmatrix}
\end{equation}
\end{La}
\vskip 0.2cm

The claim follows from Lemma \ref{LemmaTilde} along with
Eq. \eqref{eq:TildeD}.
Note that in particular,~ $\tilde{d}(s)$ ~may have multiple roots.
\vskip 0.2cm

\noindent
In the sequel, the prescribed points~ $x_1,~\ldots~,~x_m$ ~will be referred
to as~ {\em nodes}.
\vskip 0.2cm

\noindent
Without loss of generality, we shall find it convenient to
distinguish in Eq. \eqref{eq:TildeD} between the cases $b=0$, $b\not=0$.
Furthermore to ease the distinction, we differently denote
the coefficients\begin{footnote}{Although as before,~~
$\begin{smallmatrix}
{\rm Im}(x_j)>0&~&0\not=&{\gamma}_{j+1}={\gamma}_j^*
\\~\\
{\rm Im}(x_j)=0&~&0\not=&{\gamma}_j\in\R.
\end{smallmatrix}$
}\end{footnote}
of the polynomials for $b=0$ and
$b\not=0$, i.e.
\[
\begin{matrix}
\tilde{d}_o(s)&:={\tilde{d}(s)}_{|_{b=0}}&=
\sum\limits_{j=1}^m{\gamma}_j{\phi}_j(s)
\\
\tilde{d}_1(s)&:={\tilde{d}(s)}_{|_{b\not=0}}&=
b\eta(s)+\sum\limits_{j=1}^mc_j{\phi}_j(s).
\end{matrix}
\]
We now can state the following.

\begin{Tm}\label{Tm:Denom}
For distinct nodes~ $x_1,~\ldots~,~ x_m\in\C$, ~closed under complex
conjugation, let~ $\eta(s)$ ~and~
\mbox{$\phi_1(s),~\ldots~,~\phi_m(s)$} ~be as in Eqs.
\eqref{eq:eta} and \eqref{eq:phi}, respectively. Let also~
${\gamma}_1,~\ldots~,~{\gamma}_m~$ and $~c_1,~\ldots~,~c_m~$
along with $~b\in\R$, be all non-zero parameters
as in Eq. \eqref{eq:TildeD}.
\vskip 0.2cm

The set of all real polynomials~ $d(s)$ ~of degree, of at most $~m$,
with no roots at these nodes, can be parametrized by two families,
\begin{equation}\label{eq:D}
\begin{matrix}
\tilde{d}_o(s)&=&\sum\limits_{j=1}^m{\gamma}_j{\phi}_j(s)&~&~&~
&~&
{\rm deg}\left(\tilde{d}_o(s)\right)&=&m-1
\\~\\
\tilde{d}_1(s)&=&b\eta(s)+\sum\limits_{j=1}^mc_j{\phi}_j(s)&~&~&~&~&
{\rm deg}\left(\tilde{d}_1(s)\right)&=&m.
\end{matrix}
\end{equation}
\end{Tm}
\vskip 0.2cm

\subsection{Step 2:
(Not necessarily positive real)
interpolating functions with prescribed denominator.}
\label{SubSec:Interp}

\noindent
With each of the denominators in Eq. \eqref{eq:D}, $d_o(s)$ and $d_1(s)$,
we here match a numerator, denoted by $\nu_o(s)$ and $\nu_1(s)$,
respectively, to obtain a rational (not necessarily positive real)
interpolating function (from $x_j$ to $y_j$).

\begin{Tm}\label{Tm:NotNecessarilyPosInterpolating}
Let the interpolation data in Eq. \eqref{eq:interp} be
closed under complex conjugation\begin{footnote}{Namely,
if ${\rm Im}(x)_j>0$ then $x_{j+1}=x_j^*$ and
$y_{j+1}=y_j^*$.}\end{footnote}, where the nodes
$x_1,~\ldots~,~x_m$ are distinct and in $\C_l$. Let also
the (non-zero) coefficients ${\gamma}_1,~\ldots~,~{\gamma}_m$
and $c_1,~\ldots~,~c_m$ be as in Eq. \eqref{eq:D}.
\vskip 0.2cm

\noindent
Construct the polynomials
(where $\phi_j(s)$ as in Eq. \eqref{eq:phi})
\begin{equation}\label{eq:Nu}
{\nu}_o(s)=\sum\limits_{j=1}^my_j\gamma_j\phi_j(s)
\quad\quad
{\nu}_1(s)=\sum\limits_{j=1}^my_jc_j\phi_j(s)~.
\end{equation}
Then, the rational functions,
\begin{equation}\label{eq:DefP}
\begin{matrix}
\tilde{p}_o(s):=\frac{\nu_o(s)}{\tilde{d}_o(s)}
=\frac{\sum\limits_{j=1}^my_j{\gamma}_j\phi_j(s)}
{\sum\limits_{j=1}^m{\gamma}_j\phi_j(s)}
&~&
\tilde{p}_1(s):=\frac{\nu_1(s)}{\tilde{d}_1(s)}
=\frac{\sum\limits_{j=1}^my_jc_j\phi_j(s)}
{b\eta(s)+\sum\limits_{j=1}^mc_j\phi_j(s)}~,
\end{matrix}
\end{equation}
(with $0\not=b\in\R$) interpolate between $x_j$ and $y_j$.
\end{Tm}
\vskip 0.2cm

\noindent
This result follows directly from the definition
of $\phi_j(s)$ in Eq. \eqref{eq:phi}.
\vskip 0.2cm

\noindent
We next construct additional
interpolating rational functions of degree of at most $~m$.

\begin{La}\label{La:TildeFunctionR}
Let $\eta(s)$, $d_k(s)$, $\nu_k(s)$, and $p_k(s)$, (with $k=0, 1$)
from Eqs. \eqref{eq:eta}, \eqref{eq:D}, \eqref{eq:Nu}, and
\eqref{eq:DefP} respectively.\\
The set of all real rational functions of degree $~m$, vanishing at the nodes
is given by,
\begin{equation}\label{eq:TildeDelta}
\begin{matrix}
\tilde{\scriptstyle\Delta}_o(s)&=\left(\frac{\tilde{d}_o(s)}{\eta(s)}\right)^{-1}&
=\left({\sum\limits_{j=1}^m\frac{{\gamma}_j}{s-x_j}}\right)^{-1}
\\~\\
\tilde{\scriptstyle\Delta}_1(s)&=
\left(\frac{\tilde{d}_1(s)}{\eta(s)}\right)^{-1}&=
\left(b+\sum\limits_{j=1}^m
{\scriptstyle
\frac{c_j}{s-x_j}
}
\right)^{-1}
\end{matrix}
\end{equation}
with $~{\gamma}_1,~\ldots~,~{\gamma}_m$, $~b$,
$c_1, ~\ldots~,~c_m$ all non-zero.
\vskip 0.2cm

\noindent
Using $\tilde{\scriptstyle\Delta}_k(s)$,
define the rational functions,
\begin{equation}\label{eq:ParametricInterpolatingFuctionTildeF0}
\begin{matrix}
\tilde{f}_o(s)&:=&\tilde{p}_o(s)+{\scriptstyle r_o\tilde{\Delta}_o(s)}
=\frac{{\nu}_o(s)}{\tilde{d}_o(s)}+{\scriptstyle r_o}\frac{\eta(s)}{\tilde{d}_o(s)}
\\~\\~&=&
\frac{
{\sum\limits_{j=1}^my_j{\gamma}_j\phi_j(s)}+{\scriptstyle r_o}\eta(s)}
{\sum\limits_{j=1}^m{\gamma}_j\phi_j(s)}
\end{matrix}\quad\quad{\scriptstyle
r_o\in\R~~{\rm parameter}},
\end{equation}
and
\begin{equation}\label{eq:ParametricInterpolatingFuctionTildeF1}
\begin{matrix}
\tilde{f}_1(s)&:=&\tilde{p}_1(s)+{\scriptstyle r_1\tilde{\Delta}_1(s)}
=\frac{{\nu}_1(s)}{\tilde{d}_1(s)}+{\scriptstyle r_1}\frac{\eta(s)}{\tilde{d}_1(s)}
\\~\\~&=&
\frac{\sum\limits_{j=1}^my_jc_j\phi_j(s)+{\scriptstyle r_1}\eta(s)}{b\eta(s)+
\sum\limits_{j=1}^mc_j\phi_j(s)}
\end{matrix}\quad\quad{\scriptstyle
r_1\in\R~~{\rm parameter}}.
\end{equation}
Then the following is true.

\begin{itemize}
\item[(i)~~~~]{}For arbitrary $~r_o, r_1\in\R~$ the functions $~\tilde{f}_o(s)~$
and $~\tilde{f}_1(s)~$ in Eqs.
\eqref{eq:ParametricInterpolatingFuctionTildeF0}
\eqref{eq:ParametricInterpolatingFuctionTildeF1}
are of degree of at most $~m$.
\vskip 0.2cm

\item[(ii)~~]{}For arbitrary $~r_o, r_1\in\R~$ the functions $~\tilde{f}_o(s)~$
and $~\tilde{f}_1(s)~$ in Eqs.
\eqref{eq:ParametricInterpolatingFuctionTildeF0}
\eqref{eq:ParametricInterpolatingFuctionTildeF1}
interpolate from $x_1,~\ldots~,~x_m$ to $y_1,~\ldots~,~y_m$.
\end{itemize}
\end{La}

{\bf Proof :}\quad
For convenience, throughout the proof, we omit the dependence on $k=0, 1$
and simply write $\tilde{\scriptstyle\Delta}(s)$,
$\nu(s)$, $d(s)$ and $f(s)$.
\vskip 0.2cm

\noindent
$(i)$ Recall that by construction (see Eqs. \eqref{eq:eta}, \eqref{eq:D}
and Theorem \ref{Tm:Denom}) the polynomials
$\eta(s)$ and $d(s)$ are relatively prime. Recall also
(see Theorem \ref{Tm:NotNecessarilyPosInterpolating}) that
\[
m={\rm degree}\left(\eta
\right)\geq
{\rm degree}\left(\nu
\right)=m-1.
\]
We thus have the following,
\[
\begin{matrix}
m&=&{\rm degree}\left(\tilde{\scriptstyle\Delta}\right)
\\~&=&
{\rm degree}\left({\scriptstyle\frac{\eta}{d}}\right)
\\~&=&
{\rm degree}\left({\scriptstyle d}\right)+
{\rm degree}\left({\scriptstyle\eta}\right)
\\~&=&
{\rm degree}\left({\scriptstyle d}\right)+
\max\left(
{\rm degree}\left({\scriptstyle\eta}\right),
{\rm degree}\left({\scriptstyle\nu}\right)
\right)
\\~&\geq&
{\rm degree}\left({\scriptstyle d}\right)+
{\rm degree}\left({\scriptstyle\nu+r\eta}\right)
\\~&\geq&
{\rm degree}\left(
{\scriptstyle\frac{\nu+r\eta}{d}}\right)
\\~&=&
{\rm degree}\left(
\tilde{\scriptstyle f}\right).
\end{matrix}
\]
$(ii)$ This is immediate from Theorem \ref{Tm:NotNecessarilyPosInterpolating}
along with the definitions of $\tilde{\scriptstyle\Delta}(s)$ and of
$\tilde{f}(s)$, see Eqs. \eqref{eq:TildeDelta},
\eqref{eq:ParametricInterpolatingFuctionTildeF0}
and \eqref{eq:ParametricInterpolatingFuctionTildeF1}, respectively.
\qed
\vskip 0.2cm

\noindent
Following Eq. \eqref{eq:interp}, we assume hereafter
that the data set is closed under complex conjugation and that
the distinct nodes are in the open left half plane, i.e.
\[
x_1,~\ldots~,~x_m\in\C_l~.
\]
Thus, some of the resulting interpolating functions in Eqs.
\eqref{eq:ParametricInterpolatingFuctionTildeF0}
and \eqref{eq:ParametricInterpolatingFuctionTildeF1} may be positive
real while other are not. Either way, the rest of this section is
depvoted to extracting positive real functions out of them. Specifically,
we shall devise a scheme of easily constructing large subsets of
positive real interpolating functions. This simplicity, based on the
structure of interpolating functions (see Corollary
\ref{Cy:ConvexInterpolation} below), comes on the expense of guranteeing
finding ~{\em all}~  positive real interpolating functions.

\subsection{Step 3:
All positive real rational functions of degree
$~{\mathbf m}$,
with prescribed denominator, and vanishing at the nodes.}
\label{SubSec:RationalFunctions}
\vskip 0.2cm

\noindent
We first recall well-known facts, which are fundamental to our construction.
\vskip 0.2cm

\begin{Tm}\label{Tm:Convexity}
The following is true.
\begin{itemize}
\item[(i)~~]{}For prescribed data, the family of interpolating rational
functions is convex.

\item[(ii)~]{}The set of rational positive real functions forms a convex cone.
\end{itemize}
\end{Tm}
\vskip 0.2cm

As a non-empty intersection of convex sets, is convex, one
can conclude the following.

\begin{Cy}\label{Cy:ConvexInterpolation}
For prescribed data set, whenever not empty, the family of rational
positive real interpolating functions, is convex.
\end{Cy}
\vskip 0.2cm

\noindent
We now resort to the following, see e.g.
\cite[Section 4.3]{CL},
\cite[Definition 6.4]{Khalil2000}.

\begin{Dn}\label{Dn:StrictPosDef}
{\rm
A rational function $f(s)$ will be called~
{\em Strictly Positive Real}
if $~f(s-\epsilon)~$ is positive real, for some~
$\epsilon>0$.
}
\end{Dn}
\vskip 0.2cm

\noindent
The following well known properties will be useful
in the sequel.

\begin{Tm}\label{Tm:SPR}
(i) If a rational function $~f(s)~$ is strictly positive real, then
\[
{\rm Re}\left(f(s)\right)>0
\quad\quad\quad\forall s\in\overline{\C}_r~.
\]
(ii) The set of positive real functions forms a Convex Invertible
Cone.
\end{Tm}
\vskip 0.2cm

\noindent
Item (i) follows from Definition \ref{Dn:StrictPosDef}
and for item (ii) see \cite[Proposition 4.1.1]{CL}.
\vskip 0.2cm

\noindent
Whenever the functions in Eq. \eqref{eq:TildeDelta} are strictly
positive real, the tilde will be omitted and they will be denoted by~
$~{\scriptstyle\Delta}_o(s)$ and $~{\scriptstyle\Delta}_1(s)$.
This is addressed next.

\begin{Ex}\label{Ex:PosRealCoefficients}
{\rm
We next illustrate the fact that for any set of nodes in $\C_l$, the open
left-half plane, one can choose the coefficients
$~{\gamma}_1,~\ldots~,~{\gamma}_m~$ and
$~c_1,~\ldots~,~c_m$,
so that in Eq. \eqref{eq:TildeDelta} one obtains~
{\em strictly positive real}
functions, ${\scriptstyle\Delta}_o(s)$ and ${\scriptstyle\Delta}_1(s)$.
\vskip 0.2cm

\noindent
As the reasoning is identical, we show it only for
$c_1,~\ldots~,~c_m$.
\vskip 0.2cm

\noindent
Specifically, if $x_j\in\R_-$ then $\frac{c_j}{s-x_j}$ is strictly
positive real for all $c_j\in\R_+~$.\\
If $~x_j\in\{\C_l\smallsetminus\R_-\}~$ (and from
Eq. \eqref{eq:ComplexConjugation} \mbox{$x_{j+1}=x_j^*$})
then taking \mbox{$c_{j+1}=c_j^*$} yields,
\begin{equation}\label{eq:DegTwoPosDefEx}
{\scriptstyle\frac{c_j}{s-x_j}}+{\scriptstyle\frac{c_j^*}{s-x_j^*}}
=
{\scriptstyle 2{\rm Re}(c_j)}\left(s+\frac{
{\scriptstyle
\left(
{\rm Re}(c_j)
\left(-{\rm Re}(x_j)\right)
+
{\rm Im}(c_j){\rm Im}(x_j)
\right)
}
s
+
{\scriptstyle
{\rm Re}(c_j)
|x_j|^2}
}
{
{\scriptstyle
{\rm Re}(c_j)
}
s+
{\scriptstyle
\left(
{\rm Re}(c_j)
\left(-{\rm Re}(x_j)\right)
-
{\rm Im}(c_j){\rm Im}(x_j)
\right)
}
}
\right)^{-1}.
\end{equation}
Thus, choosing $c_j$ so that,
\[
{\scriptstyle{\rm Re}(c_j)}>
{\scriptstyle\frac{|{\rm Im}(c_j){\rm Im}(x_j)|}{-{\rm Re}(x_j)}}\geq 0,
\]
is sufficient to guarantee that the function in Eq. \eqref{eq:DegTwoPosDefEx}
is strictly positive real.
}
\qed
\end{Ex}
\vskip 0.2cm

\noindent
Motivated by the above example we can state the following.
\vskip 0.2cm

\begin{Tm}\label{Tm:Delta}
Out of the polynomials $~\eta(s)$, $~\tilde{d}_o(s)~$ and $~\tilde{d}_1(s)~$
in Eqs. \eqref{eq:eta}, \eqref{eq:D} respectively, one can choose the
coefficients $~c_1~,~\ldots~,~c_m$, $b~$ and
$~{\gamma}_1~,~\ldots~,~{\gamma}_m~$ to construct all
monic\begin{footnote}{As in Eqs.
\eqref{eq:ParametricInterpolatingFuctionTildeF0}
\eqref{eq:ParametricInterpolatingFuctionTildeF1}
$~\tilde{\scriptstyle\Delta}(s)$ is scaled by $r$,
without loss of generality one can take $d(s)$
to be moinc.}\end{footnote} polynomials,
\begin{equation}\label{eq:Denominator}
\begin{matrix}
{d}_o(s)=&\sum\limits_{j=1}^m{\gamma}_j{\phi}_j(s)&
\sum\limits_{j=1}^m{\gamma}_j=1&
{\rm deg}\left({d}_o(s)\right)=&m-1
\\~\\
{d}_1(s)=&\eta(s)+\sum\limits_{j=1}^mc_j{\phi}_j(s)&~&
{\rm deg}\left({d}_1(s)\right)=&m,
\end{matrix}
\end{equation}
so that the rational functions from Eq. \eqref{eq:TildeDelta} are
strictly positive real, i.e.
\begin{equation}\label{eq:Delta}
\begin{matrix}
{\scriptstyle\Delta}_o(s)=
\left(\frac{{d}_o(s)}{\eta(s)}\right)^{-1}=
\left(
{\scriptstyle\sum\limits_{j=1}^m\frac{{\gamma}_j}{s-x_j}}
\right)^{-1}
&\sum\limits_{j=1}^m{\gamma}_j=1
\\~\\
{\scriptstyle\Delta}_1(s)=\left(\frac{{d}_1(s)}{\eta(s)}\right)^{-1}=
\left(1+\sum\limits_{j=1}^m
{\scriptstyle
\frac{c_j}{s-x_j}
}
\right)^{-1}.
&c_1{\cdots}c_m\neq 0.
\end{matrix}
\end{equation}
Moreover:\\
The coefficients of $~{d}_o(s)$: ${\gamma}_1~,~\ldots~,~{\gamma}_m$ form
a convex set (excluding ${\gamma}_j=0$) within a
hyperplane in $\C^m$.\\
The coefficients of $~{d}_1(s)$: $c_1,~\ldots~,~c_m$ form a convex subset of
$\C^m$ (excluding $~m~$ hyperplanes $c_j=0$).
This set is positively unbounded in the sense that if
in Eqs. \eqref{eq:D}, \eqref{eq:Delta}
\[
c_1,~\ldots~,~c_m
\]
is an admissible set of parameters, then so is\begin{footnote}{
provided that, to preserve complex conjugation,
$c_{j+1}=c_j^*$ implies $\delta_{j+1}=\delta_j~$.}\end{footnote}
\[
c_1+\delta_1,~\ldots~,~c_m+\delta_m
\quad\quad\quad\quad
\delta_j\geq 0\quad j=1,~\ldots~,~m,
\]
Furthermore, the sets ${\gamma}_1~,~\ldots~,~{\gamma}_m$ and
$c_1,~\ldots~,~c_m$ can be
parameterized by a convex subset of $\R^m$ (excluding the axes).
\end{Tm}

{\bf Proof}\quad
From Example \ref{Ex:PosRealCoefficients} it follows that for arbitrary set of
nodes, this family of ${\scriptstyle\Delta}(s)$ functions is not empty.
\vskip 0.2cm

\noindent
To simplify establishing structural properties, we begin by
ignoring the condition that neither \mbox{${\gamma}_1,~\ldots~,~{\gamma}_m$}
nor $c_1,~\ldots~,~c_m$ vanish.
\vskip 0.2cm

\indent
By Corollary \ref{Cy:ConvexInterpolation}
the set of interpolating positive real rational functions is convex.
\vskip 0.2cm

\noindent
Indeed, if $c_1~,~\ldots~,~c_m$ and $\hat{c}_1~,~\ldots~,~\hat{c}_m$ are
two admissible sets in Eqs. \eqref{eq:D}, \eqref{eq:Delta} then
so is ${\scriptstyle\theta}{c_1}+(1-{\scriptstyle\theta})\hat{c}_1~
,~\ldots~,~{\scriptstyle\theta}{c_m}+(1-{\scriptstyle\theta})\hat{c}_m$,
for all ${\scriptstyle\theta}\in[0,~1]$.
\vskip 0.2cm

\noindent
The fact that $\sum\limits_{j=1}^m{\gamma}_j=1$ forms a hyper-plane
in $\C^m$ is straightforward. Next, to show that the set
$c_1,~\ldots~,~c_m$ is positively unbounded, one can resort again to
the construction in Example \ref{Ex:PosRealCoefficients}.
\vskip 0.2cm

\noindent
Recall however that the restriction that $\prod\limits_{j=1}^mc_j\not=0$,
implies that the set of coefficients $~c_1,~\ldots~,~c_m~$ forms
an almost convex cone, as it excludes $~m~$ hyper-planes, $~c_j=0$.
\vskip 0.2cm

\noindent
Real coefficients: Recall that  in \eqref{eq:D} the coefficients $c_j$ are
real or come in complex conjugate pairs. Specifically if
there are $q$ coefficients in the upper half plane and
$m-2q$ are real, they are described by a point in $\R^m$.
\qed
\vskip 0.2cm

\noindent
The above analysis suggests that in the coefficient space, it is enough
to find the {\em boundary} of the (almost convex) sets of admissible
${\gamma}_1,~\ldots~,~{\gamma}_m~$ and $~c_1,~\ldots~,~c_m$.
\vskip 0.2cm

\noindent
In the next step, we combine Theorems \ref{Tm:NotNecessarilyPosInterpolating}
and \ref{Tm:Delta}
to construct positive real interpolating functions of degree of at most $m$.
\vskip 0.2cm

\subsection{Step 4: Positive real
interpolating functions.}
\label{SubSec:PosInterp}

\noindent
To extract positive real functions, out of the set of interpolating
functions $\tilde{f}(s)$ in Eqs.
\eqref{eq:ParametricInterpolatingFuctionTildeF0} and
\eqref{eq:ParametricInterpolatingFuctionTildeF1},
we focus on those whose deniminator is given by Theorem \ref{Tm:Delta}.
This is formalized next.

\begin{La}\label{La:FunctionR}
Let $\eta(s)$, $d_k(s)$, $\nu_k(s)$, $p_k(s)$ and $~{\scriptstyle\Delta}_k(s)$,
(with $k=0, 1$) from Eqs. \eqref{eq:eta}, \eqref{eq:D}, \eqref{eq:Nu},
\eqref{eq:DefP} and \eqref{eq:TildeDelta} respectively.\\
Define the rational functions,
\begin{equation}\label{eq:ParametricInterpolatingFuctionF0}
\begin{matrix}
{f}_o(s)&:=&p_o(s)
+
{\scriptstyle
r_o
{\Delta}_o(s)}
=
\frac{{\nu}_o(s)}{d_o(s)}
+
{\scriptstyle
r_o
}
\frac{\eta(s)}{d_o(s)}
\\~\\~&=&
\frac{
{\sum\limits_{j=1}^my_j{\gamma}_j\phi_j(s)}
+{\scriptstyle r_o}\eta(s)}
{\sum\limits_{j=1}^m{\gamma}_j\phi_j(s)}
\end{matrix}\quad\quad{\scriptstyle
r_o\in\R~~{\rm parameter}},
\end{equation}
and
\begin{equation}\label{eq:ParametricInterpolatingFuctionF1}
\begin{matrix}
{f}_1(s)&:=&p_1(s)+{\scriptstyle r_1{\Delta}_1(s)}=
\frac{{\nu}_1(s)}{d_1(s)}+{\scriptstyle r_1}\frac{\eta(s)}{d_1(s)}
\\~\\~&=&
\frac{\sum\limits_{j=1}^my_jc_j\phi_j(s)+{\scriptstyle r_1}\eta(s)}
{b\eta(s)+\sum\limits_{j=1}^mc_j\phi_j(s)}
\end{matrix}\quad\quad{\scriptstyle
r_1\in\R~~{\rm parameter}}.
\end{equation}
Then, $~f_k(s)$ (and
$p_k(s)$) are interpolating function with $~{\scriptstyle\Delta}_k(s)$
strictly positive real, vanishing at the nodes (all sharing the same denominator).
\vskip 0.2cm

\noindent
Furthermore, the quantities,
\[
\begin{matrix}
-\inf\limits_{s\in{\C}_r}~\frac
{
{\rm Re}~p_o(s)
}
{
{\rm Re}~{\scriptstyle\Delta}_o(s)
}
&~&~&~&~&~&
-\inf\limits_{s\in{\C}_r}~\frac
{
{\rm Re}~p_1(s)
}
{
{\rm Re}~{\scriptstyle\Delta}_1(s)
}
\end{matrix}
\]
are well defined.
\end{La}

{\bf Proof :}\quad
For convenience, throughout the proof, we omit the dependence on $k=0, 1$
and simply write $\nu(s)$, $d(s)$ and $f(s)$.
\vskip 0.2cm

\noindent
The construction in Theorem \ref{Tm:Delta} guarantees that
in Eq. \eqref{eq:Delta}
\[
{\scriptstyle\Delta}(s)={\frac{\eta(s)}{d(s)}}
\]
is~ {\em strictly positive real} and thus by item (i) of
Theorem \ref{Tm:SPR}
\[
{{\rm Re}\frac{\eta(s)}{d(s)}}>0\quad\quad\quad
\forall{ s}\in{\C}_r~.
\]
Recalling that (i) the numerator $~\eta(s)$ vanishes only at $~m~$
points in ${\C}_l$, see Eq. \eqref{eq:eta} and (ii)
in addition ${\rm degree}(\eta)\geq{\rm degree}(d)$
see Eq. \eqref{eq:D}, in fact
\[
{{\rm Re}\frac{\eta(s)}{d(s)}}\geq\delta>0
\quad\quad\quad\forall{ s}\in{\C}_r~.
\]
Next, exploiting again the fact that $~{\scriptstyle\Delta}(s)~$ is
strictly positive real, see Theorem \ref{Tm:SPR}, implies that
$~d(s)~$ does not vanish in $\overline{\C}_r$. We can thus conclude
that
\[
\inf\limits_{s\in{\C}_r}~\frac
{{\rm Re}~p(s)}
{{\rm Re}~{\scriptstyle\Delta}(s)}
\]
is well defined, so the claim is established.
\qed
\vskip 0.2cm

\noindent
So far we have described interpolating rational functions $~f(s)~$
of degree of at most $~m$ whose denominator is so that
$~{\scriptstyle\Delta}(s)~$ is strictly positive real.
To proceed with the construction, the idea is very simple, see
Eq. \eqref{eq:KeyIdea}~:\\
With the same $~\eta(s)$, $~d(s)~$ and $~\nu(s)$ construct
the rational functions
\[
f_a(s):=
{\frac{\nu(s)}{d(s)}}
+
{
r_a
\frac{\eta(s)}{d(s)}}
\quad\&\quad
f_b(s):=
{\frac{\nu(s)}{d(s)}}
+
{ r_b\frac{\eta(s)}{d(s)}}
\]
where $r_a$ and $r_b$ are real parameters.
On the one hand, from Lemma \ref{La:FunctionR} it follows that $f_a(s)$
and $f_b(s)$ interpolate between with same data. On the other hand,
Theorems
\ref{Tm:SPR}
and
\ref{Tm:Delta}
imply that,
\[
r_a>r_b\quad\Longrightarrow\quad
{\rm Re}\left(f_a(s)\right)
>
{\rm Re}\left(f_b(s)\right)
\quad s\in\C_r~.
\]
Next, recall that by Eq. \eqref{eq:DefPos} $~f(s)~$ is positive real
whenever,
\[
{\rm Re}\left(f(s)\right)\geq 0
\quad\quad\quad\quad\forall s\in\C_r~.
\]
Thus, one can formally define
$\forall s\in\C_r$,
\begin{equation}\label{eq:DefUnderlineR}
\begin{matrix}
\underline{r}_o:=
{\rm arg}~\min\limits_{r_o\in\R}~{\rm Re}\left(f_o(s)\right)
={\rm arg}~\min\limits_{r_o\in\R}~{\rm Re}\left(p_o(s)
+r_o{\scriptstyle\Delta}_o(s)\right)\geq 0
\\~\\
\underline{r}_1:={\rm arg}~\min\limits_{r_1\in\R}~{\rm Re}\left(f_1(s)\right)
={\rm arg}~\min\limits_{r_1\in\R}~{\rm Re}\left(p_1(s)+
r_1{\scriptstyle\Delta}_1(s)\right)\geq 0.
\end{matrix}
\end{equation}
We next combine the above definition of $\underline{r}$ along with
Lemma \ref{La:FunctionR}.
\vskip 0.2cm

\begin{Pn}\label{Pn:UnderlineR}
Let the rational function $~f_k(s)~$ and the scalars
$~\underline{r}_k$ (with $k=0, 1$) be as in Eqs.
\eqref{eq:ParametricInterpolatingFuctionF0}, \eqref{eq:ParametricInterpolatingFuctionF1}
and \eqref{eq:DefUnderlineR}, respectively.\quad Then,
\[
\begin{matrix}
\underline{r}_o:=-\inf\limits_{s\in{\C}_r}~
\frac{{\rm Re}~p_o(s)}{{\rm Re}~{\scriptstyle\Delta}_o(s)}
&~&~&~&~&~&
\underline{r}_1:=-\inf\limits_{s\in{\C}_r}
\frac{{\rm Re}~p_1(s)}{{\rm Re}~{\scriptstyle\Delta}_1(s)}
\end{matrix}
\]
and $~f_k(s)~$ is positive real if and only if $~r_k\geq\underline{r}_k~$.
\end{Pn}

\noindent
{\bf Proof :}\quad For simplicity, we omit both the dependence on $~s~$
and the subscript $~k$. Using Eqs.
\eqref{eq:ParametricInterpolatingFuctionF0}
\eqref{eq:ParametricInterpolatingFuctionF1} note that
\[
{\rm Re}\left(f\right)={\rm Re}
\left(
{
\frac{\nu}{d}+r\frac{\eta}{d}}\right)
={\rm Re}\left(
{
\frac{\nu}{d}}\right)+r{\rm Re}\left(
{
\frac{\eta}{d}}\right).
\]
Now, $~f~$ is positive real if and only if
\[
{\rm Re}\left(f\right)\geq 0
\quad\quad\quad
\forall s\in{\C}_r~.
\]
Namely,
\[
r{\rm Re}\left({\frac{\eta}{d}}\right)\geq
-{\rm Re}\left({\frac{\nu}{d}}\right)
\quad\quad\quad
\forall s\in{\C}_r~,
\]
in turn, using the fact that ${\frac{\eta}{d}}$ is strictly
positive real, see Theorem \ref{Tm:Delta}, this means that
\[
r\geq
{\frac{-{\rm Re}~
\frac{\nu}{d}~
}
{~~~{\rm Re}~
\frac{\eta}{d}~
}}
\quad\quad\quad
\forall s\in{\C}_r~.
\]
Hence, one can conclude that $~f~$ in Eqs.
\eqref{eq:ParametricInterpolatingFuctionF0},
\eqref{eq:ParametricInterpolatingFuctionF1}
is
positive real, if and only if,
\[
r\geq
\sup\limits_{s\in{\C}_r}~
{\frac{-{\rm Re}~
\frac{\nu}{d}~
}
{~~~{\rm Re}~
\frac{\eta}{d}~
}}
=
-\inf\limits_{s\in{\C}_r}~
{\frac{{\rm Re}~
\frac{\nu}{d}~
}
{{\rm Re}~
\frac{\eta}{d}~
}}~,
\]
and by Eq. \eqref{eq:DefUnderlineR}, the proof is complete.
\qed
\vskip 0.2cm

\noindent
Noting that ${\rm deg}(\eta_k)=m-1$, for $k=0, 1$, while
${\rm deg}(\psi)=m$ together with the fact that
${\scriptstyle\Delta}_k$ is strictly positive real,
guarantees the following.

\begin{Ob}\label{Ob:rGeq0}
In Proposition \ref{Pn:UnderlineR},
\[
\underline{r}_o\geq 0\quad\quad\quad\underline{r}_1\geq 0.
\]
\end{Ob}

\noindent
Note that from Eqs.  \eqref{eq:ParametricInterpolatingFuctionF0}
and \eqref{eq:ParametricInterpolatingFuctionF1} it follows that
for $~r_k>0$, with $k=0, 1$ whenever there is no pole-zero cancelation,
the degree of the numerator of $~f_o(s)~$ or of $~f_1(s)~$ is $~m$.
Thus, all positive real interpolating functions $~f(s)~$ we have
constructed are of degree at most $~m$, but  {\em under the restriction
that the degree of the numerator is greater or equal to the degree of
the denominator}. In the next section we address the complementary
case where the degree of the
denominator is greater or equal to the degree of the numerator.

\subsection{Step 5: Additional positive real
interpolating functions.}
\label{SubSec:AllPosInterp}

Taking the original data, if one considers a function, say $g(s)$,
interpolating from \mbox{$x_1,~\ldots~,~x_m$} to\begin{footnote}{Assuming
$y_j\not=0$}\end{footnote}
\mbox{$\frac{1}{y_1}~,~\ldots~,~\frac{1}{y_m}$}, then $~\frac{1}{g(s)}~$
solves the original problem, where we have relied on the fact that the inverse
of a positive real function, is positive real, see item (ii) of
Theorem \ref{Tm:SPR}. Here are the details.
\vskip 0.2cm

\noindent
We follow the previous steps (while adding hat to the respective
functions) and first
mimic Theorem \ref{Tm:NotNecessarilyPosInterpolating}.

\begin{Tm}\label{Tm:HatNu}
Let the interpolation data\begin{footnote}{Assuming
$~y_j\not=0$.}\end{footnote} be as in Eq. \eqref{eq:interp}, the
(non-zero) denominator coefficients ${\gamma}_1,~\ldots~,~{\gamma}_m~$
and $~c_1,~\ldots~,~c_m~$ be as in Eq. \eqref{eq:D}.
Construct the polynomials
(where $\phi_j(s)$ as in Eq. \eqref{eq:phi})
\[
\hat{\nu}_o(s)=\sum\limits_{j=1}^m\frac{{\gamma}_j}{y_j}\phi_j(s)
\quad\quad\hat{\nu}_1(s)=\sum\limits_{j=1}^m
\frac{\hat{c}_j}{y_j}\phi_j(s).
\]
Then, the rational functions,
\[
\hat{p}_o(s):=\frac{\hat{\nu}_o(s)}{d_o(s)}
=
\frac
{\sum\limits_{j=1}^m\frac{{\gamma}_j}{y_j}\phi_j(s)}
{\sum\limits_{j=1}^m{\gamma}_j\phi_j(s)}
\quad\quad
\quad
\hat{p}_1(s):=\frac{\hat{\nu}_1(s)}{d_1(s)}
=
\frac
{\sum\limits_{j=1}^m\frac{c_j}{y_j}\phi_j(s)}
{\eta(s)+\sum\limits_{j=1}^mc_j\phi_j(s)}
~,
\]
interpolate between $~x_j~$ and $~\frac{1}{y_j}~$ with
$~j=1,~\ldots~,~m$.
\end{Tm}
\vskip 0.2cm

\noindent
Note that indeed all the parameters are as before.
\vskip 0.2cm

\noindent
We next mimic Lemma \ref{La:FunctionR} and
construct the rational functions,
\begin{equation}\label{eq:ParametricInverseInterpolationHatF0}
\begin{matrix}
\hat{f}_o(s)&:=&
\left(\hat{p}_o(s)+
{
\hat{r}_o{\scriptstyle\Delta}_o(s)}
\right)^{-1}
&~&{\hat{r}_o\in\R~~{\rm parameter}}
\\~\\~&=&
\left(
\frac{\hat{\nu}_o(s)}{d_o(s)}+
{
\hat{r}_o}\frac{\eta(s)}{d_o(s)}
\right)^{-1}
&=&
\left(
\frac{
\sum\limits_{j=1}^m\frac{{\gamma}_j}{y_j}\phi_j(s)
+\hat{r}_o\eta(s)}
{
\sum\limits_{j=1}^m{\gamma}_j\phi_j(s)
}
\right)^{-1}
\end{matrix}
\end{equation}
and
\begin{equation}\label{eq:ParametricInverseInterpolationHatF1}
\begin{matrix}
\hat{f}_1(s)&:=&
\left(\hat{p}_1(s)
+
{
\hat{r}_1
{\scriptstyle\Delta}_1(s)
}
\right)^{-1}
&~&{\hat{r}_1\in\R~~{\rm parameter}}
\\~\\~&=&
\left(
\frac{\hat{\nu}_1(s)}{d_1(s)}+{\hat{r}_1}\frac{\eta(s)}{d_1(s)}
\right)^{-1}
&=&
\left(
\frac{
\sum\limits_{j=1}^m\frac{c_j}{y_j}\phi_j(s)
+{\hat{r}_1}\eta(s)
}
{\eta(s)+\sum\limits_{j=1}^mc_j\phi_j(s)}
\right)^{-1}
\end{matrix}
\end{equation}
where $\eta(s)$, $d_k(s)$, ${\scriptstyle\Delta}_k(s)$ and $p_k(s)$, with
$k=0, 1$ are ~{\em as before}, see Eqs.  \eqref{eq:eta}, \eqref{eq:D}
\eqref{eq:Delta} and \eqref{eq:DefP}, respectively.
\vskip 0.2cm

\noindent
As before, the rational function $~\hat{f}_k(s)~$ (with $k=0, 1$)
see Eqs.
\eqref{eq:ParametricInverseInterpolationHatF0},
and
\eqref{eq:ParametricInverseInterpolationHatF1}
interpolate
from $x_1,~\ldots~,~x_m$ to $y_1,~\ldots~,~y_m$
for all $~\hat{r}_k\in\R$. Out of this family, we next extract the
positive real subset. To this end,
we introduce the following notation,
\begin{equation}\label{eq:DefUnderlineHatR}
\begin{smallmatrix}
\underline{\hat{r}}_o:=
{\rm arg}~\min\limits_{\hat{r}_o\in\R}~{\rm Re}\left(\hat{f}_o(s)\right)
&=
{\rm arg}~\min\limits_{\hat{r}_o\in\R}~{\rm Re}\left(
\frac{\hat{\nu}_o(s)}{d_o(s)}+
\hat{r}_o\frac{\eta(s)}{d_o(s)}\right)^{-1}\geq 0
\\~\\
\underline{\hat{r}}_1:=
{\rm arg}~\min\limits_{\hat{r}_1\in\R}~{\rm Re}\left(\hat{f}_1(s)\right)
&=
{\rm arg}~\min\limits_{\hat{r}_1\in\R}~{\rm Re}\left(
\frac{\hat{\nu}_1(s)}{d_1(s)}+
\hat{r}_1\frac{\eta(s)}{d_1(s)}\right)^{-1}\geq 0
\end{smallmatrix}
\quad\quad\forall s\in\C_r~.
\end{equation}
\noindent
By using item (ii) of Theorem \ref{Tm:SPR},
we can next adapt Proposition \ref{Pn:UnderlineR}
to guarantee that the sought interpolating functions are indeed positive real.

\begin{Pn}\label{Pn:HatUnderlineR}
Let the rational function $~\hat{f}_k(s)~$ and the scalars
$~\underline{\hat{r}}_k(s)$ (with $k=0, 1$) be as in Eqs.
\eqref{eq:ParametricInverseInterpolationHatF0},
\eqref{eq:ParametricInverseInterpolationHatF1}
and
\eqref{eq:DefUnderlineHatR}, respectively.\quad Then,
\[
\begin{matrix}
\underline{\hat{r}}_o:=-\inf\limits_{s\in{\C}_r}~
\frac{{\rm Re}~\hat{p}_o(s)
}
{{\rm Re}~
{\scriptstyle\Delta}_o(s)
}
&~&~&~&~&~&
\underline{\hat{r}}_1:=-\inf\limits_{s\in{\C}_r}
\frac{{\rm Re}~\hat{p}_1(s)
}
{{\rm Re}~{\scriptstyle\Delta}_1(s)
}
\end{matrix}
\]
and $~\hat{f}_k(s)$ is positive real if and only if
$~\hat{r}_k\geq\underline{\hat{r}}_k~$.
\end{Pn}
\vskip 0.2cm

\noindent
We have shown that $~\hat{f}_o(s)~$ and $~\hat{f}_1(s)~$ are
positive real interpolating functions of degree at most $~m$,
where the degree of the numerator is larger or equal to the
degree of the denominator.
\vskip 0.2cm

\noindent
Simlar to the reasoning at end of Step 4, one can conclude
the following.

\begin{Ob}\label{Ob:HatrGeq0}
In Proposition \ref{Pn:HatUnderlineR},
\[
\underline{\hat{r}}_o\geq 0\quad\quad\quad\underline{\hat{r}}_1\geq 0.
\]
\end{Ob}

\section{Examples and Concluding remarks}
\setcounter{equation}{0}

The above recipe is illustrated through simple examples.
\vskip 0.2cm

\noindent
{\bf A.}\quad We start by illustrating the role of $~f_o(s)~$
vs. $~\hat{f}_o(s)$ in Eqs.
\eqref{eq:ParametricInterpolatingFuctionF0}
and \eqref{eq:ParametricInverseInterpolationHatF0}, respectively
to obtain interpolating functions having at $s=\infty$
either pole or zero.
\vskip 0.2cm

\noindent
{\bf (i)}\quad
Find a minimal degree positive real function $~f(s)~$ mapping
$x_1,~\ldots~,~x_m\in\C_l$ to
\mbox{$y_1=x_1,~\ldots~,~y_m=x_m$}.
Clearly the sought solution is
\[
f(s)=s.
\]
We now follow the above recipe and substitute in Eq.
\eqref{eq:ParametricInterpolatingFuctionF0}
\[
\begin{split}
f_o(s)&=\dfrac{{\sum\limits_{j=1}^my_j{\gamma}_j\phi_j(s)}
+r_o\eta(s)}
{\sum\limits_{j=1}^m{\gamma}_j\phi_j(s)}\\
&=\frac{\sum\limits_{j=1}^mx_j{\gamma}_j\phi_j(s)+r_o\eta(s)}
{\sum\limits_{j=1}^m{\gamma}_j\phi_j(s)}\quad\hspace{4.1cm}\mbox{\text{\rm for $y_j=x_j$}}\\
&=\frac
{\sum\limits_{j=1}^m
x_j{\gamma}_j\phi_j(s)+
\frac{r_o}{m}
\sum\limits_{j=1}^m
(s-x_j)
\phi_j(s)}
{\sum\limits_{j=1}^m{\gamma}_j\phi_j(s)}
\\
&=
\frac
{
\sum\limits_{j=1}^m
\left(s+x_j(m\gamma_j-1)\right)
\phi_j(s)}
{m\sum\limits_{j=1}^m{\gamma}_j\phi_j(s)}\quad\hspace{3.3cm}\mbox{\text{\rm for $r_o=1$}}
\\
&=s \quad\hspace{8cm}\mbox{\text{\rm for $\gamma_j\equiv\frac{1}{m}~.$}}
\end{split}
\]
\vskip 0.2cm

{\bf (ii)}
Find a minimal degree positive real function $~f(s)~$ mapping
\mbox{$x_1,~\ldots~,~x_m\in\C_l$} to
\mbox{$y_1=\frac{1}{x_1},~\ldots~,~y_m=\frac{1}{x_m}$}.
Clearly the sought solution is
\[
f(s)={\dfrac{1}{s}}~.
\]
We now follow the above recipe and substitute in Eq.
\eqref{eq:ParametricInverseInterpolationHatF0}
\[
\begin{split}
\hat{f}_o(s)&=\left(\frac{\sum\limits_{j=1}^m\frac{{\gamma}_j}{y_j}\phi_j(s)
+\hat{r}_o\eta(s)}{\sum\limits_{j=1}^m{\gamma}_j\phi_j(s)}
\right)^{-1}
\\&
=
\left(
\frac
{\sum\limits_{j=1}^mx_j{\gamma}_j\phi_j(s)
+\hat{r}_o\eta(s)}
{\sum\limits_{j=1}^m{\gamma}_j\phi_j(s)}
\right)^{-1}
\quad\hspace{3.1cm}\mbox{\text{\rm for $y_j=\frac{1}{x_j}$}}\\
&
=
\left(
\frac
{
\sum\limits_{j=1}^m
\left(s+x_j(m\gamma_j-1)\right)
\phi_j(s)}
{m\sum\limits_{j=1}^m{\gamma}_j\phi_j(s)}
\right)^{-1}
\quad\hspace{2.2cm}\mbox{\text{\rm for $\hat{r}_o=1$}}\\
&=\frac{1}{s}
\quad\hspace{8cm}\mbox{\text{\rm for $\gamma_j\equiv\frac{1}{m}~.$}}
\end{split}
\]
\vskip 0.2cm

\noindent
{\bf B.}\quad
Parametrize all positive real rational functions, of degree of
at most two, so that
\[
f(-1)=y_1\quad\quad\quad\quad f(-3)=y_2~,
\]
where $y_1, y_2\in\R$ are arbitrary.
\vskip 0.2cm

\noindent
First for reference, a direct computation reveals that all
rational functions, of degree of at most one, are given by
\begin{equation}\label{eq:ExReafDegreeOne}
f(s)=\frac{
{\scriptstyle\left(a(3y_2-y_1)+b(y_1-y_2)\right)}s+
{\scriptstyle 3a(y_2-y_1)+b(3y_1-y_2)}
}
{{\scriptstyle 2}({\scriptstyle a}s+{\scriptstyle b})}~.
\end{equation}
These functions are positive real whenever,
\begin{equation}\label{eq:ExReafDegreeOneConstraints}
\begin{matrix}
a&\geq 0\\~\\
b&\geq 0\\~\\
a(3y_2-y_1)+b(y_1-y_2)&\geq 0\\~\\
3a(y_2-y_1)+b(3y_1-y_2)&\geq 0.
\end{matrix}
\end{equation}
The conditions in Eq. \eqref{eq:ExReafDegreeOneConstraints} may be
satisfied for all $~y_1,~y_2\in\R$ ~{\em umless},
$~0>y_1=y_2$.
\vskip 0.2cm

\noindent
This implies that for $y_1=y_2\geq 0$ there is a zero degree positive
real interpolating function, see item (iii) below. For $~0>y_1=y_2$,
the positive real interpolating functions are of degree of at least two,
see item (vi) below. In all other cases, there exist positive real
interpolating functions of degree one and above.
\vskip 0.2cm

\noindent
We now follow the recipe from the previous section.

From Step 1
\[
\eta(s)=(s+1)(s+3)=s^2+4s+3
\]
and
\[
\phi_1(s)=s+3\quad\quad\phi_2(s)=s+1.
\]
From Step 2, and using Eq. \eqref{eq:D} yields
\[
d_o(s)=s+\gamma\quad\quad\quad
{\gamma\in[0,~4]\smallsetminus\{1, 3\}},
\]
and
\[
{ d_1(s)}=
{
s^2+s(4+c_1+c_2)+3+3c_1+c_2}
\]
where $c_1$ and $c_2$ are such that
$ 1+{\frac{c_1}{s+1}}+
{\frac{c_2}{s+3}}$ is strictly positive real.
For $d_1(s)$ the set of admissible parameters is convex and
positively unbounded\begin{footnote}{From Theorem \ref{Tm:Delta}
it follows that in particular it contains the whole first quadrant
of the $\{c_2,c_1\}$ plane}\end{footnote}
(excluding the axes $c_1=0$ and $c_2=0$), it is given by
\begin{equation}\label{eq:c1c2}
c_2>\left\{\begin{matrix}-3(c_1+1)&~&\frac{1}{8}\geq c_1\\~\\
-\frac{1}{3}(\sqrt{c_1}+{ 2\sqrt{2}})^2&~&c_1\geq\frac{1}{8}~.
\end{matrix}\right.
\quad\quad\quad c_1c_2\neq 0.
\end{equation}
\begin{figure}[ht!]
\begin{tikzpicture}[scale=1.2,cap=round]
\pgfplotsset{compat=1.10}
\pgfplotsset{every axis/.append style={
                   axis x line=middle,    
                   axis y line=middle,    
                   axis line style={->,color=black, thick}, 
                   xlabel={$c_1$},          
                   ylabel={$c_2$},          
              }}
     \begin{axis}[
             xmin=-5.1,xmax=35.3,
             ymin=-24.1,ymax=15.8,
             grid=both,
             ]
\addplot [domain=-5.0:0.125,samples=100,color=gray, very thin, name path=A]({x}, {12});
\addplot [domain=-5.0:0.125,samples=100,color=blue, thick, name path=B]({x},{-3*x-3});
\addplot [gray] fill between[of=A and B];

\addplot [domain=0.125:32,samples=1000,color=gray, very thin, name path=C]{12};
\addplot [domain=0.125:32,samples=1000,color=blue, thick, name path=D]({x},{-(x+4*(2*x)^0.5+8)/3});
\addplot [gray] fill between[of=C and D];
     \end{axis}
\end{tikzpicture}
\caption{\mbox{$c_1$, $c_2$ for $\frac{\eta}{d_1}$ strictly positive real,
Eq. \eqref{eq:c1c2}.
}}
\label{Figure:parameters}
\end{figure}

From Step 3
\[
\begin{matrix}
\frac{{\nu}_o(s)}{d_o(s)}&=&
{\frac
{
y_1(\gamma-1)
}{2}}
+
{\frac
{
y_2(3-\gamma)
}{2}}
+
{\frac
{(\gamma-3)(\gamma-1)(y_2-y_1)}
{2(s+\gamma)}}
&~&
{\gamma\in[0,~4]\smallsetminus\{1, 3\}}
\\~\\
\frac{{\nu}_1(s)}{d_1(s)}&=&\frac{(c_1y_1+c_2y_2)s
+3c_1y_1+c_2y_2}{s^2+s(4+c_1+c_2)+3+3c_1+c_2}&~&c_1,~c_2
\quad{\rm from~Eq.~\eqref{eq:c1c2}}.
\end{matrix}
\]
Now from Step 4
\begin{equation}\label{eq:ExNo}
f_o(s)=r_o(s+4-{\scriptstyle\gamma})+
{\scriptstyle\frac{y_1(\gamma-1)}{2}}+
{\scriptstyle\frac{y_2(3-\gamma)}{2}}
+(r_o+{\scriptstyle\frac{y_2-y_1}{2}})
\frac{\scriptstyle(\gamma-3)(\gamma-1)}{s+{\scriptstyle\gamma}}
\end{equation}
with $\gamma\in[0,~4]\smallsetminus\{1, 3\}$, and with ${ c_1,~c_2}$
from Eq. \eqref{eq:c1c2},
\begin{equation}\label{eq:ExN1}
f_1(s)=
\frac
{
{
r_1}s^2+
{
(4r_1+c_1y_1+c_2y_2)}s+
{
3r_1+3c_1y_1+c_2y_2
}
}
{
s^2+
{ (4+c_1+c_2)}s+
{ 3+3c_1+c_2}
}~.
\end{equation}
One can verify that taking $~r_o$, $~r_1~$ ``sufficiently large"
renders $f_o(s)$, $~f_1(s)~$ positive real.
\vskip 0.2cm

\noindent
Next, from Step 5, assuming $~{ y_1y_2}\not=0$
and $\gamma\in[0,~4]\smallsetminus\{1, 3\}$,
\begin{equation}\label{eq:ExHatFo}
\hat{f}_o(s)=\frac
{\scriptstyle
{ 2y_1y_2}(s+\gamma)
}
{\scriptstyle
{ 2y_1y_2\hat{r}_o}s^2+
{
(8y_1y_2\hat{r}_o+(3-\gamma)y_1+(\gamma-1)y_2)}s
+{ 6y_1y_2\hat{r}_o+(3-\gamma)y_1+3(\gamma-1)y_2}}
\end{equation}
and with ${ c_1,~c_2}$
from Eq. \eqref{eq:c1c2},
\begin{equation}\label{eq:ExHatF1}
\hat{f}_1(s)=
\frac
{
s^2+
{
(4+c_1+c_2)
}
s+
{
3+3c_1+c_2
}
}
{
{
\left(\frac{c_1}{y_1}+\frac{c_2}{y_2}\right)}s+
{
\frac{3c_1}{y_1}+\frac{c_2}{y_2}
}
+
{
\hat{r}_1
(s+1)(s+3)
}
}
\end{equation}
Again, taking $~\hat{r}_o$, $~\hat{r}_1$ ``sufficiently large" renders
$~\hat{f}_o(s)$, $~\hat{f}_1(s)~$ positive real.
\vskip 0.2cm

\noindent
Here are five particular cases.
\vskip 0.2cm

\noindent
{\bf (i)}\quad
Recall that in the Introduction we pointed out that if $y_1, y_2\in\R_-$,
see Eq. \eqref{eq:OddNevPick}, one can still try to resort to the
classical Nevanlinna-Pick interpolation, seeking positive real odd
fuctions so that
\[
f(-1)=y_1\quad\quad f(1)=-y_1 \quad\quad f(-3)=y_2 \quad\quad f(3)=-y_2~.
\]
Now, the solvability condition in Eq.
\eqref{eq:NevPick} reads,
\begin{equation}\label{eq:CondExOddNevPick}
{\frac{y_2}{y_1}}\in\left[
{\scriptstyle\frac{1}{3}}~,~3\right],
\end{equation}
and the resulting positive real odd interpolating functions
(of degree of at most two) are
\[
\begin{matrix}
g_a(s)&=&
\frac{8y_1y_2s}{(y_2-3y_1)s^2+3(y_1-3y_2)}
\\~\\
g_b(s)&=&
\frac{\left(y_1-3y_2\right)s^2+3\left(y_2-3y_1\right)}{8s}~.
\end{matrix}
\]
We now show, that these positive real odd functions, are special cases
of the above recipe:\\
Indeed, assuming the condition in Eq.
\eqref{eq:CondExOddNevPick} is staisfied, from Eq. \eqref{eq:ExN1}
\[
{f_1(s)}_{|_{c_1=\frac{4y_2}{3y_1-y_2}~~c_2=\frac{12y_1}{y_2-3y_1}~~r_1=0}}
=
g_a(s),
\]
and from Eq. \eqref{eq:ExHatF1}
\[
{\hat{f}_1(s)}_{|_{c_1=\frac{4y_1}{3y_2-y_1}~~c_2=\frac{12y_2}{y_1-3y_2}~~\hat{r}_1=0}}
=g_b(s).
\]
\vskip 0.2cm

\noindent
To further emphasize that our approach is different, in the four
following special cases (ii), (iii) and (v), the condition in Eq.
\eqref{eq:CondExOddNevPick} is not satisfied, so the classical
Nevanlinna-Pick interpolation is not applicable.
\vskip 0.2cm

\noindent
{\bf (ii)}\quad Take the special case where $y_1=1$ and $y_2=3$.\\
Clearly, $f(s)=-s$ is a real, anti-positive, minimal degree,
interpolating function. We next seek minimal
degree positive real interpolating functions.
\vskip 0.2cm

\noindent
Substituting these image points in $~f_o(s)~$ in Eq.
\eqref{eq:ExNo} yields the following positive real
interpolating functions,
\[
f_o(s)={ r_o}s+
{({r_o}+1)}\left(4-{\scriptstyle\gamma}
+\frac{\scriptstyle (\gamma-3)(\gamma-1)}
{s+{\scriptstyle\gamma}}\right)
\quad\quad
{\gamma\in[0,~4]\smallsetminus\{1, 3\}}.
\]
To guarantee minimal degree, further substitute
$~{\underline{r}_o}=0$, to obtain interpolating
functions with zero at infinity,
\[
f_o(s)=4-{\scriptstyle\gamma}
+\frac{\scriptstyle(\gamma-3)(\gamma-1)}
{s+{\scriptstyle\gamma}}
\quad\quad
{\gamma\in[0,~4]\smallsetminus\{1, 3\}}.
\]
Comparing with Eqs. \eqref{eq:ExReafDegreeOne} and
\eqref{eq:ExReafDegreeOneConstraints} reveals that in this case
our recipe yields all minimal degree (equals one) positive real
interpolating functions.
\vskip 0.2cm

\noindent
Similarly for $~\hat{f}_o(s)~$ in Eq. \eqref{eq:ExHatFo}
\[
\hat{f}_o(s)={\frac{s+{\scriptstyle\gamma}}
{
{\left(\frac{\gamma-1}{2}+\frac{3-\gamma}{6}\right)}s+
{
\frac{3(\gamma-1)}{2}+\frac{3-\gamma}{6}
}
+
{\hat{r}_o}(s+1)(s+3)}
}
\quad\quad\quad
\]
with $\gamma\in[0,~4]\smallsetminus\{1, 3\}$.
As before, to single out interpolating functions of degree one,
we focus on cases where $~\underline{\hat{r}}_o=0$. However,
then to guarantee positive realness, the range of the parameter
$~{\gamma}~$ needs to be further restricted, i.e.
\[
\hat{f}_o(s)=3\left({\scriptstyle\gamma}+
\frac{\scriptstyle(3-\gamma)(\gamma-1)}
{s+{\scriptstyle\gamma}}\right)^{-1}
\quad\quad
{\gamma\in\left[{\scriptstyle\frac{3}{4}},~4\right]\smallsetminus\{1, 3\}}.
\]
Here, at infinity, the interpolating function has neither pole nor zero.
\vskip 0.2cm

\noindent
Finally note that comparison with Eqs. \eqref{eq:ExReafDegreeOne}
and \eqref{eq:ExReafDegreeOneConstraints} reveals that in this case,
the recipe produced ~{\em all}~ interpolating functions of degree
one.
\vskip 0.2cm

\noindent
{\bf (iii)}\quad Take the special case where $y_1=y_2\geq 0$.\\
One can substitute in Eq. \eqref{eq:ExNo} $~\underline{r}_o=0~$
to obtain the minimal (=zero) degree interpolating
function \mbox{$f_o(s)\equiv y_1$}.\\
Similarly, one can substitute in
Eq. \eqref{eq:ExHatFo} $~\underline{\hat{r}}_o=0~$
to obtainn the minimal (=zero) degree interpolating
function \mbox{$\hat{f}_o(s)\equiv y_1$}.
\vskip 0.2cm

\noindent
{\bf (iv)}\quad
Take the special case where $0>y_1=y_2~$.\\
Recall that from Eqs. \eqref{eq:ExReafDegreeOne} and
\eqref{eq:ExReafDegreeOneConstraints} we know that there are no positive
real interpolating function of degree less than two.
\vskip 0.2cm

\noindent
To obtain interpolating functions use the recipe and
substitute in Eq. \eqref{eq:ExNo} to obtain,
\[
f_o(s)=y_1+r_o\frac{(s+3)(s+1)}{s+\gamma}
\quad\quad\underline{r}_o=\left\{
\begin{smallmatrix}
-y_1\frac{\gamma}{3}&~&~&\gamma\in(1,~3)\\~\\
-y_1\frac{1}{4-\gamma}&~&~&\gamma\in\{[0,~1)\cup(3,~4)\}.
\end{smallmatrix}
\right.
\]
Note that $\underline{r}_o$ turns to be unbounded, as $\gamma$
approaches 4.
\vskip 0.2cm

\noindent
{\bf (v)}\quad
Take the special case where $y_1=2$, $y_2=0$.\\
As before, substituting these image points
in Eqs.  \eqref{eq:ExReafDegreeOne} and \eqref{eq:ExReafDegreeOneConstraints}
(with $\frac{b}{a}={\scriptstyle\gamma}$)
reveals that {\em all}~ minimal degree (equals one) positive real interpolating
functions are of the form\begin{footnote}{Substituting in Eq. \eqref{eq:ExNo}
$r_o=0$, yields the subset of the interpolating functions
in Eq. \eqref{eq:ExZeroImageDegOne}, where
$4\geq{\scriptstyle\gamma}$}\end{footnote}
\begin{equation}\label{eq:ExZeroImageDegOne}
f(s)=({\scriptstyle\gamma}-1)\frac{s+3}{s+{\scriptstyle\gamma}}
\quad\quad\quad 3\neq{\scriptstyle\gamma}>1.
\end{equation}
Next, address the case where the interpolating function is so that
the degree of the denominator is strictly larger then the the degree
of the numerator.
Now, recall that since the set of image points contains zero,
Step 5 of the recipe cannot be used. Nevertheless, all required
interpolating functions are obtained.
\vskip 0.2cm

\noindent
We start with a straightforward considerations: Since at $x=-1$, the
numerator is non-zero, but it vanishes at $x=-3$, it must be (at least)
of degree one. Thus, the denominator is (at least) of degree two.
Indeed, to obtain all minimal degree interpolating functions of the
required nature, substitute in Eq. \eqref{eq:ExN1}
\[
{f_1(s)}_{|_{r_1=0}}=
\frac{2c_1(s+3)}{(s+1)(s+3)+c_1(s+3)+c_2(s+1)}~,
\]
where adapting Eq.\eqref{eq:c1c2},
\[
0\neq c_2>\left\{\begin{matrix}-3(c_1+1)&~&c_1\in(0, \frac{1}{8}]\\~\\
-\frac{1}{3}(\sqrt{c_1}+{ 2\sqrt{2}})^2&~&c_1\geq\frac{1}{8}~.
\end{matrix}\right.
\]
\vskip 0.2cm

\noindent
{\bf C}.~
In the previous item the interpolation nodes were real.
We here illustrate the fact that the recipe is identical
for the non-real case, assuming the interpolation nodes
are closed under complex conjugation.
\vskip 0.2cm

\noindent
Assume that the interpolation nodes are
\mbox{$x_1=-\gamma+i\delta$} and
\mbox{$x_2=-\gamma-i\delta$} where $\gamma>0$ and $0\not=\delta\in\R$.
Hence,
\[
\eta(s)=(s+\gamma)^2+\delta^2
\]
We now construct the denominator polynomials.
\vskip 0.2cm

\noindent
Following Theorem \ref{Tm:Delta} a degree one numerator polynomial
$~d_o(s)~$ is given by the condition that the following rational
function is strictly positive real,
\[
\begin{matrix}
\dfrac{\eta(s)}{d_o(s)}
&=&\left(\dfrac{\scriptstyle\frac{1}{2}+i\beta}{s+\gamma+i\delta}+
\dfrac{\scriptstyle\frac{1}{2}+i\beta}{s+\gamma+i\delta}\right)^{-1}\\~\\~&=&
\dfrac{(s+\gamma)^2+{\delta}^2}{s+\gamma+2\beta\delta}\\~\\~&=&
s+\gamma-2\beta\delta+\dfrac{\scriptstyle\delta^2(1+4\beta^2)}{
s+\gamma+2\beta\delta}~,
\end{matrix}
\]
namely,
\[
\gamma>2|\beta\delta|.
\]
Hence one arrives at the following parametrization,
\[
d_o(s)=s+2\gamma(1-\theta)
\quad\quad\quad\quad\theta\in[0,~1).
\]
\qed
\vskip 0.2cm

{\bf Concluding remarks}
\vskip 0.2cm

\noindent
1.\quad As already pointed out in Corollary \ref{Cy:ConvexInterpolation},
for arbitrary prescribed data set in $\C$, the family of all
positive real interpolating functions is
convex (whenever not empty).
\vskip 0.2cm

\noindent
In contrast,
the set of rational functions of a degree of at most $m$ is a cone,
but highly non-convex. In fact, the degree of a sum of two rational
functions is higher than the degree of each of the summands,
unless one of the denominators divides the other.
\vskip 0.2cm

\noindent
When the interpolation nodes are in $\C_l$, the open left half plane,
we here introduce an easy-to-compute parametrization of positive real
interpolating functions as a subset of $\R^{2m+3}$, see item 2 for
details.
\vskip 0.2cm

\noindent
2.\quad For arbitrary interpolating data set in Eq. \eqref{eq:interp},
closed under complex conjugation, with nodes in $\C_l$, a large
subset of positive real interpolating functions of degree of at most
$m$ may be conveniently parametrized a union of convex subsets within
$\R^{2m+3}$.
\vskip 0.2cm

\noindent
Indeed the coefficients in Eq. \eqref{eq:Denominator}
are so that $c_1,~\ldots~,~c_m$ form a positively unbounded convex subset
of $\R^m$, which in particular contains $\R_+^m$,
excluding the axes (see e.g. Figure \ref{Figure:parameters}). Next,
${\gamma}_1,~\ldots~,~{\gamma}_m$ form a hyper-plane in $\R^{m-1}$.
Finally, each of the four parameters $\underline{r}_o$, $\underline{r}_1$,
$\underline{\hat{r}}_o$,
$\underline{\hat{r}}_1$,
lies in $\overline{\R}_+$.
\vskip 0.2cm

\noindent
3.\quad Step 4 of the recipe relies on the fact that positive real
rational functions form a convex cone and that the set of interpolating
functions is convex. Steps 3 and 5 rely  on the fact that the set
of positive real rational functions is closed under inversion.
\vskip 0.2cm

\noindent
4.\quad The parametrization through
$f_o(s)$, $f_1(s)$, $\hat{f}_o(s)$, $\hat{f}_1(s)$
is motivated by ~{\em simplicity}. It is neither minimal, as
the same interpolation function may be obtained in
more than one way, see e.g. Example B(iii), nor is it comprehensive,
as some of the minimal degree interpolating functions
may be missing, see e.g. Example B(v).
\vskip 0.2cm

\noindent
5.\quad While the parametrization through $f_o(s)$, $f_1(s)$,
$\hat{f}_o(s)$, $\hat{f}_1(s)$ is convenient, focusing on~
{\em minimal degree} interpolating functions involves ``fine
tuning" of the parameters $~{\gamma}_1,~\ldots~,~{\gamma}_m$,
$~c_1, ~\ldots~,~c_m$, $~r_o$, $~r_1$, $~\hat{r}_o~$ and
$~\hat{r}_1$, see Examples A, B.
\vskip 0.2cm

\begin{center}
Acknowledgement
\end{center}

The authors thank Prof. V. Bolotnikov form the Math. Dept. at
the College of William and Mary, Williamsburg, Virginia, USA
for providing useful, constructive remarks at early stage of
this work.

\end{document}